\documentstyle[12pt]{article}

\input amssym.def
\input amssym

\newtheorem{dfn}{Definition}[section]

\newtheorem{lem}[dfn]{Lemma}

\newtheorem{prop}[dfn]{Proposition}
\newtheorem{cor}[dfn]{Corollary}

\newtheorem{gdt}[dfn]{Geometric Decomposition Theorem}

\def\proof{\par\medskip\noindent{\it Proof: }}
\def\qed{{\unskip\nobreak\hfil
        \penalty50\hskip1em\hbox{}\nobreak\hfil
        $\square$\parfillskip=0pt\finalhyphendemerits=0 \par}}

\def\ga{\gamma}
\def\Ga{\Gamma}

\def\si{\sigma}
\def\Si{\Sigma}

\def\R{{\Bbb R}}

\def\MIN{\mathop{\hbox{\rm MIN}}}

\begin{document}

\title{A Geometric Characteristic Splitting in all Dimensions}
\author{by Bernhard Leeb\thanks{Partially supported by SFB 256 (Bonn) and IHES} 
and Peter Scott\thanks{Partially supported by NSF grant DMS 9306240}}
\date{Preprint - May 20, 1996}
\maketitle

\abstract{\noindent We prove the existence of a geometric characteristic
submanifold for non-positively curved ma\-ni\-folds of any dimension greater
than or equal to three. In dimension three, our result is a geometric
version of the topological characteristic submanifold theorem due to Jaco,
Shalen and Johannson.} \bigskip

\section{Introduction}

In the 1970's, Jaco and Shalen \cite{JacoShalen} and Johannson \cite
{Johannson} showed that a closed orientable Haken $3$-manifold $M$ has a
canonical family of disjoint embedded incompressible tori, no two of which
are parallel, such that the complementary pieces of $M$ are either Seifert
fibre spaces or are atoroidal. They defined the characteristic submanifold $%
V(M)$ of $M$ to be essentially the union of the Seifert manifold pieces of $%
M $. Further, they showed that any essential map of the torus into $M$ is
homotopic into $V(M)$. Johannson called this last property the Enclosing
Property. For brevity, we will refer to these results as the JSJ results.

In this paper, we show that if $M$ is a closed manifold of dimension three
or more, and if $M$ has a Riemannian metric of non-positive curvature, then
either the metric on $M$ is flat or there is a precisely analogous
decomposition of $M$ along codimension one submanifolds. Further these
submanifolds are totally geodesic in $M$ and are flat in the metric induced
from $M$. Note that in dimension three, a flat manifold must be a Seifert
fibre space, so that, in particular, our arguments give a new proof of the
JSJ results for the special case when $M$ is assumed to have a metric of
non-positive curvature. In dimension four or more, a flat manifold need not
be a Seifert manifold, see the example near the end of section \ref{prelim},
so this case really is different in higher dimensions. We also prove that
essentially the same results hold if $M$ is non-orientable and if the
boundary of $M$ is non-empty, on the assumption that the boundary is convex.

At the time when Jaco and Shalen and Johannson proved their results, the
methods seemed very special to dimension three and no one even asked whether
this result had any generalization to higher dimensions. Several years later
in 1990, Kropholler \cite{Kropholler} published an algebraic analogue of
their results. He considered Poincar\'e duality groups of dimension three ($%
PD3$-groups). The fundamental group of any closed aspherical $3$-manifold is
automatically a $PD3$-group, but it is not known whether the converse holds.
Note, however, that $PD2$-groups are known to be fundamental groups of
closed surfaces, so it is not unreasonable to hope that every $PD3$-group is
the fundamental group of an aspherical closed $3$-manifold. There is a
natural analogue, in the context of $PD3$-groups, of an embedded
incompressible surface. One considers a $PD3$-group $G$ and a subgroup $H$
such that $H$ is a $PD2$-group and $G$ splits over $H$, i.e. $G$ can be
expressed as $A*_H$, or as $A*_HB$ with $A\neq H\neq B$. Kropholler showed
that the natural analogue of the JSJ splitting result holds for $PD3$%
-groups. Surprisingly, he also showed that his algebraic result had a
generalization to $PDn$-groups in all dimensions greater than three. This
raised the question of whether the topological results of Jaco, Shalen and
Johannson also generalized to higher dimensions. As Haken manifolds are
aspherical and Kropholler's results correspond to results about aspherical
manifolds, it seems possible that the JSJ results might generalise to
aspherical manifolds but not to all manifolds. Note that a Riemannian
manifold of non-positive curvature with convex boundary is aspherical, so
our results show that the JSJ results generalize to aspherical manifolds in
the special case of non-positive curvature.

Jaco and Shalen \cite{JacoShalen} and Johannson \cite{Johannson} also
considered non-closed manifolds and defined a characteristic submanifold $%
V(M)$ for any orientable Haken $3$-manifold $M$ with incompressible
boundary. They showed that such a manifold $M$ has a canonical family of
disjoint properly embedded incompressible tori and annuli, no two of which
are parallel. They defined the characteristic submanifold $V(M)$ of $M$ to
be essentially the union of the Seifert manifold pieces of $M$ together with
some pieces which are homeomorphic to $I$-bundles. Further, they showed that
any essential map of the torus or annulus into $M$ is homotopic into $V(M)$.
In this paper, we give analogous results for compact manifolds with
non-empty boundary in any dimension greater than or equal to three, but we
assume that $M$ has a metric of non-positive curvature and that the
boundary of $M$ is totally geodesic.

The results in this paper were proved by the authors independently in 1992.
At about the same time, Sela \cite{Sela} announced some algebraic results
which are closely related to all the preceding discussion. Sela's results
were for negatively curved groups and were the precise analogue of the JSJ
results for the case of $3$-manifolds with no incompressible tori. This is
because a negatively curved group cannot have a subgroup isomorphic to ${\bf %
Z\times Z}$. The topological picture is of a $3$-manifold with a canonical
family of disjoint embedded annuli, and Sela's picture is of a group which
splits over several different infinite cyclic subgroups. Sela's Enclosing
Property is the analogue of the JSJ Enclosing Property for embedded annuli
only. More recently, Rips and Sela \cite{RipsSela} have announced a
generalization of Sela's results to cover all finitely presented groups.
Thus again we are left with the question of whether the results in this
paper can be generalised to manifolds which need not have non-positive
curvature.

As we pointed out earlier, a Riemannian manifold $M$ of non-positive
sectional curvature with convex boundary is aspherical and so its homotopy
type is determined by its fundamental group. It is well-known that various 
algebraic properties of $\pi_1(M)$ have strong implications for the geometry 
of $M$.
The most basic example, due to Gromoll-Wolf \cite{GromollWolf} 
and Lawson-Yau \cite{LawsonYau}, 
is that an abelian subgroup of $\pi_1(M)$ 
is carried by a totally-geodesically immersed flat torus. 
In this paper, 
we obtain information about the geometric structure 
of non-positively curved manifolds
from the intersection pattern of the closed flat 
totally-geodesic hypersurfaces.
Our main results are geometric versions of the topological
decomposition theorem in dimension three due to Jaco, Shalen \cite
{JacoShalen} and Johannson \cite{Johannson}.

\medskip
\noindent
{\bf Geometric Decomposition Theorem in Dimension Three.} {\em Let }$M$ {\em %
be a compact connected non-positively curved 3-manifold which has convex
boundary. Then either }$M$ {\em is closed and has a flat metric, or }$M${\em %
\ can be canonically decomposed along finitely many
totally-geo\-de\-si\-cal\-ly embedded flat 2-tori and Klein bottles. The
resulting pieces are Seifert or atoroidal. } {\em Further any }$\pi _1${\em %
-injective map of the torus or Klein bottle into }$M${\em \ can be homotoped
to a totally geodesic flat immersion, and any such immersion must lie in one
of the Seifert pieces or be homotopic to a cover of one of the 
decomposing surfaces.}\medskip

Note that some of the decomposing surfaces may be one-sided. 
In particular, no piece in the decomposition of $M$ will be an interval bundle 
over a flat surface unless $M$ itself is an interval bundle over a flat surface. 
If $M$ is a twisted interval bundle over a flat surface $F$ 
then our construction splits $M$ along the one-sided surface $F$. 

The Seifert pieces of $M$ admit a Seifert fibration by closed geodesics and
they are rigid in the sense that they split locally as a Riemannian product,
the fiber being the one-dimensional factor. Note that if $M$ is flat, it is
also Seifert fibered in the three dimensional case. 
The proof of our theorem readily applies to all dimensions. 
See the end of section \ref
{prelim} for the definitions.

\medskip
\noindent
{\bf Geometric Decomposition Theorem.} {\em Let }$M$ {\em be a compact
connected non-positively curved manifold of dimension }$n\geq 3${\em , which
has convex boundary. Then either }$M$ {\em is closed and has a flat metric,
or }$M${\em \ can be canonically decomposed along finitely many
totally-geo\-de\-si\-cal\-ly embedded flat closed submanifolds of
codimension one. The resulting pieces are Seifert fibered or codimension-one
atoroidal. } {\em Further any essential manp of a closed flat }$(n-1)${\em %
-manifold into }$M${\em \ can be homotoped to a totally geodesic flat
immersion, and any such immersion must lie in one of the Seifert pieces 
or be homotopic to a cover of one of the 
decomposing hypersurfaces.}\medskip

As in the three-dimensional case, 
some of the decomposing hypersurfaces may be one-sided. 
There is also a more general version of this result which corresponds to the
full JSJ splitting of an orientable compact $3$-manifold along annuli as
well as tori. We leave the statement to section \ref{general case}.

This paper is organized as follows: In section~\ref{submanifolds} we prove
that there is an upper bound for the number of mutually non-parallel,
disjoint, totally-geodesically embedded, closed hypersurfaces in a compact
non-positively curved manifold $M$ with convex boundary. 
In section~\ref
{rigid geometry}, we 
study the pattern ${\cal S}$ of totally-geodesically immersed, flat, closed
hypersurfaces in $M$. 
We show that intersecting hypersurfaces span a geometric Seifert
submanifold. The desired decomposition of $M$ is obtained by cutting along
hypersurfaces in ${\cal S}$ which are {\em isolated} in the sense that they do
not intersect any other surface in ${\cal S}$. 
In section \ref{general case}, we discuss how to
prove the most general version of our results.

\section{Preliminaries}
\label{prelim}

\noindent
{\bf Non-positive curvature.} 
We start by recalling a few well-known facts 
from the geometry of nonpositively curved manifolds, 
for more details the reader may consult e.g.\ \cite{CheegerEbin}. 
In this paper, we consider smooth Riemannian manifolds 
of non-positive sectional curvature.
We will always assume that they are complete as metric spaces 
and that their boundaries are convex, 
i.e.\ each
geodesic touching the boundary must already be contained in the boundary. 
A simply-connected manifold $X$ of this kind
has the fundamental property that its 
distance function $d:X\times X\to {\Bbb R}$
is {\em convex}, 
that is,  
for any two geodesics $c_1,c_2:[a,b]\to X$ the function $%
t\to d(c_1(t),c_2(t))$ is convex. In particular, the distance $d(\cdot ,C)$
from a convex subset $C\subset X$ is a convex function. The convexity of $d$
implies that any two points in $X$ can be connected by a {\em unique}
geodesic. This has strong topological implications: $X$ is contractible so that 
any manifold $M$ covered by $X$ is aspherical, 
i.e.\ $M$ is a $K(\pi ,1)$-space.

A smooth submanifold $Y\subseteq X$, possibly with boundary, is called {\em %
totally-geodesic} if each geodesic in $X$ tangent to an interior point of $Y$
belongs locally to $Y$. We call $Y$ {\em (geodesically) complete} 
if each geodesic in $Y$
is extendable ad infinitum. If $Y_1$ and $Y_2$ are complete totally-geodesic
submanifolds of $X$ which have bounded distance from each other, then the
distance functions $d(\cdot ,Y_i)\mid _{Y_j}$ are constant by convexity and
completeness. 
This implies that the submanifolds $Y_i$ are {\em parallel},
i.e.\ there is a totally-geodesic submanifold in $X$ which splits metrically
as $Y\times [a_1,a_2]$ so that $Y_i=Y\times \{a_i\}$. 

For an isometry $\phi$ of $X$, denote by $\MIN(\phi)$ 
the set where the {\em displacement
function} $d_{\phi}:x\to d(x,\phi x)$ assumes its infimum. 
Since $d$ is convex, $%
\mathop{\hbox{\rm MIN}}(\phi)$ is a closed convex subset of $X$. 
An isometry $\phi$ is called {\em non-parabolic} or {\em semisimple} if $\MIN(\phi)$ 
is non-empty. 
In this case, $\phi$ is  {\em elliptic} if the minimum of $d_{\phi}$ equals zero
and {\em loxodromic} if it is strictly positive. 
The set of minimal displacement for a loxodromic isometry splits metrically as
\begin{equation}
\mathop{\hbox{\rm MIN}}(\phi)\cong{\Bbb R}\times Y 
\label{splitting of minimal set}
\end{equation}
where the lines ${\Bbb R}\times \{y\}$ are the $\phi$-{\em axes},
i.e.\  geodesics preserved by $\phi$, 
and $Y$ is a simply-connected manifold of nonpositive curvature 
with convex boundary. 
Isometries of 
$X$ commuting with $\phi$ preserve the splitting (\ref
{splitting
of
minimal
set}). 
It follows by induction that
any abelian subgroup $A$ of the isometry group of $X$ preserves a flat in $X$
where a {\em flat} is defined to be a convex subset 
isometric to a Euclidean space. 
More precisely, 
the intersection of minimal sets $\bigcap_{\gamma \in A}%
\mathop{\hbox{\rm MIN}}(\ga)=:\mathop{\hbox{\rm MIN}}(A)$ is non-empty and
splits metrically as 
\begin{equation}
\label{splittingofminsetofabeliabgroup}
MIN(A)\cong E\times Y
\end{equation}
where $E$ is a Euclidean space 
(possibly of dimension zero) and 
$Y$ is a simply-connected manifold of nonpositive curvature 
with convex boundary. 
The layers $E\times \{y\}$ are the minimal $A$-invariant flats and  
the induced action of $A$ on $E$ is cocompact.

Suppose that $\Ga$ is a group which acts properly-discontinuously and cocompactly
by isometries on $X$, 
such as the group of deck-transformations corresponding to a compact 
Riemannian manifold covered by $X$. 
Let $A\subset\Ga$ be an abelian subgroup 
(which is necessarily finitely generated) 
and denote by $C(A)$ its centraliser and by $N(A)$ its normaliser. 
The action of $N(A)$ on $X$ preserves $%
\mathop{\hbox{\rm MIN}}(A)$ and the splitting 
(\ref{splittingofminsetofabeliabgroup}). 

\begin{lem}
\label{centraliseractscoco}
The action of $C(A)$ on $\mathop{\hbox{\rm MIN}}(A)$ is cocompact.
\end{lem}

\medskip\noindent{\it Proof: } 
Let $(p_n)$ be a sequence of points in $%
\mathop{\hbox{\rm MIN}}(A)$. 
Since $\Gamma $ acts cocompactly on $X$,
there exist isometries $\gamma _n\in \Gamma $ so that the sequence $(\gamma
_np_n)$ is bounded. Let $a_1,\ldots ,a_r$ denote a basis of $A.$ For each
value of the index $i$, the points 
\[
\gamma _na_i\gamma _n^{-1}\cdot \gamma _np_n=\gamma _n\cdot a_ip_n, 
\]
form a bounded sequence too, because $d(a_ip_n,p_n)$ equals the minimal
displacement of the isometry $a_i$. Since the action of $\Gamma $ is
properly discontinuous, the elements $\gamma _na_i\gamma _n^{-1}$ are
contained in a finite subset of $\Gamma $. By passing to a subsequence $r$
times, we can assume that, for each $i$, $\gamma _na_i\gamma _n^{-1}$ is a
fixed element of $\Gamma $ for all values of $n$. Let $\gamma _n^{^{\prime
}} $ denote $\gamma _1^{-1}\gamma _n$ which must lie in $C(A)$. Then the
sequence $(\gamma _n^{\prime }p_n)$ is bounded, as it is obtained from the
bounded sequence $(\gamma _np_n)$ by applying $\gamma _1^{-1}$. Since $%
(p_n)\subset \mathop{\hbox{\rm MIN}}(A)$ was chosen arbitrarily, we conclude
that there is a bounded fundamental domain for the action of $C(A)$ on $%
\mathop{\hbox{\rm MIN}}(A)$.
\qed

\noindent
The following auxiliary result will be needed later:

\begin{lem}
\label{flats in products} 
Let $X_1$ and $X_2$ be simply-connected Riemannian
manifolds of non-positive curvature 
(metrically complete and with convex boundary). If $F\subset
X_1\times X_2$ is a totally-geodesically embedded flat submanifold, then the
images of $F$ under the projections $p_i:X_1\times X_2\rightarrow X_i$ on
the factors are also flat.
\end{lem}

\medskip\noindent{\it Proof: } Let $c$ and $c^{\prime }$ be geodesic
segments in $F$ so that their distance function $d(t):=d(c(t),c^{\prime
}(t)) $ is constant. Denote by $d_i$ the distance function of the projected
segments $p_i\circ c$ and $p_i\circ c^{\prime }$. Then $d^2=d_1^2+d_2^2$.
Since the $d_i$ are convex, $d^2$ can only be constant if the $d_i$ are
constant. Hence $p_i$ maps parallel segments to parallel segments and the
claim follows.
\qed

\bigskip\noindent
{\bf Topology.} We explain the notions necessary to state the topological
decomposition theorem due to Jaco, Shalen and Johannson. We work in the
smooth category. Let $M$ be a compact orientable 3-manifold, possibly with
boundary, which is {\em irreducible}, i.e.\ every embedded 2-sphere bounds
an embedded 3-ball, and has infinite fundamental group. We consider
connected, two-sided, embedded surfaces $\Sigma $ in $M$ which are not
homeomorphic to the 2-sphere. We will also require that $\Sigma $ be
properly embedded in $M$ or be embedded in the boundary of $M$. Such a
surface $\Sigma $ is called {\em incompressible} if there is no disc $D$
embedded in $M$ such that $D\cap \Sigma =\partial D$ and $\partial D$ is a
non-contractible curve in $\Sigma $. If $M$ contains a properly embedded
incompressible surface, then $M$ is a {\em Haken} manifold. The following
decomposition theorem has been proven for Haken manifolds by Jaco, Shalen 
\cite{JacoShalen} and Johannson \cite{Johannson}. The non-Haken case follows
from the fact that if a compact orientable irreducible $3$-manifold admits a 
$\pi _1$-injective map of the torus but does not admit such an embedding,
then it must be a Seifert fibre space. This result requires the work of
several authors and the proof was completed independently by Casson and
Jungreis \cite{CassonJungreis} and Gabai \cite{Gabai}.

\medskip\noindent
{\bf Topological Decomposition Theorem.} {\em A compact orientable
irreducible 3-manifold with infinite fundamental group and incompressible
boundary can be cut along finitely many disjoint incompressible 2-tori into
atoroidal and Seifert pieces, and any }$\pi _1${\em -}${\em injective}${\em %
\ map of the 2-torus into the manifold is homotopic into one of the Seifert
pieces or to a covering of one of the decomposing tori. 
Moreover, a minimal such decomposition is unique up to isotopy.}

It remains to explain the types of pieces which occur: These are compact
3-manifolds $N$ with boundary. A $3$-manifold $N$ is {\it atoroidal }if any $%
\pi _1$-injective map of the torus into $N$ is homotopic into the boundary
of $N.$ It is a {\em Seifert} manifold if it admits a Seifert fibration,
i.e.\ if it can be expressed as a disjoint union of embedded circles, the
fibres, so that the following is true: Every fibre has a neighborhood which
is isomorphic, as a fibred space, to a fibred solid torus or Klein bottle. A
fibred solid torus is a quotient of the trivially fibred product $D^2\times 
{\Bbb R}$ by a diffeomorphism $(\phi ,\tau )$ where $\phi $ is an isometry
of finite order of the unit disc $D^2$ and $\tau $ is a translation on the
real line.

We also need to define what is meant by the terms Seifert manifold and
atoroidal in higher dimensions. In dimension three, a Seifert manifold is a
Seifert bundle over a $2$-dimensional orbifold with fiber the circle. In the
context of this paper, we define a Seifert manifold $N$ of dimension $n$ to be a
Seifert bundle over a $2$-dimensional orbifold with fiber
a flat $(n-2)$-manifold. 
This means that $N$ is foliated by $(n-2)$-dimensional closed flat manifolds 
so that each leaf has a foliated neighborhood which has a finite cover 
whose induced foliation is a product $F\times D^2$.  
A manifold $M$ of dimension $n$ is {\it %
codimension-one atoroidal} if any $\pi _1$-injective map of a flat $(n-1)$%
-torus into $M$ is homotopic into the boundary of $M$.

\section{Immersed totally-geodesic submanifolds}
\label{submanifolds}

>From now on, 
$M$ will denote a compact connected Riemannian manifold $M$ of
non-positive curvature with convex boundary. We denote by $\pi :\tilde{M}\to
M$ the universal covering map and think of $\pi _1(M)=:\Gamma $ as group of
deck transformations acting on $\tilde{M}$.

Let $\phi :\Sigma \to M$ be a totally-geodesic Riemannian immersion of a
closed connected non-positively curved manifold $\Sigma $ into $M$. 
Every lifting to a
map of universal covers is a totally-geodesic embedding $\tilde \phi :\tilde
\Sigma \hookrightarrow \tilde M$ and induces an injective
homomorphism $\pi _1(\Sigma )\hookrightarrow \Gamma =\pi _1(M)$ of
fundamental groups. Different lifts yield conjugate subgroups of $\Gamma $.
Note that $\tilde\Si$ is geodesically complete. 

\subsection{Intersections}

\begin{lem}
\label{cocompactintersections}
Let $C_1,C_2\subset\tilde M$ be closed subsets 
so that the stabiliser $\Ga_i:=Stab_{\Ga}(C_i)$ acts cocompactly on $C_i$. 
Then $\Ga_1\cap\Ga_2$ acts cocompactly on $C_1\cap C_2$. 
\end{lem}
\proof
The natural map 
$(\Ga_1\cap\Ga_2)\backslash\Ga_2\to\Ga_1\backslash\Ga$ is injective
and the corresponding immersion
$(\Ga_1\cap\Ga_2)\backslash C_2\to \Ga_1\backslash\tilde M$ 
is therefore proper. 
Hence the inverse image under this immersion 
of the compact subset $\Ga_1\backslash C_1$ 
is compact.
As this inverse image equals 
$(\Ga_1\cap\Ga_2)\backslash(C_1\cap C_2)$,
the lemma follows.  
\qed

Note that the lemma holds more generally 
for properly discontinuous group actions on locally-compact topological spaces.

\begin{cor}
Let $\Si_1$ and $\Si_2$ be closed non-positively curved Riemannian manifolds
and suppose that $\phi_1:\Sigma _1\to M$ and $\phi_2:\Sigma _2\to M$
are totally-geodesic Riemannian immersions. 
Then $\phi_1(\Sigma_1)\cap\phi_2(\Sigma_2)$
is a finite union of totally-geodesically immersed 
closed non-positively curved Riemannian manifolds.
\end{cor}
\proof
The immersion $\phi_i$ lifts to an embedding of universal covers 
with image a closed convex subset $Y_i\subset\tilde M$. 
By the previous lemma,
the totally-geodesic submanifolds $\ga_1\cdot Y_1\cap\ga_2\cdot Y_2$,
$\ga_1,\ga_2\in\Ga$, have cocompact stabilisers in $\Ga$.
The corollary follows because, by compactness, 
$\phi_1(\Sigma_1)\cap\phi_2(\Sigma_2)$ is the projection of finitely many 
submanifolds $\ga_1\cdot Y_1\cap\ga_2\cdot Y_2$.
\qed

\subsection{Finiteness for disjoint non-parallel totally-geo\-de\-sic hypersurfaces}

\begin{dfn}
We call two totally-geodesic Riemannian immersions $\phi_1:\Sigma _1\to M$ 
and $\phi_2:\Sigma _2\to M$ 
of closed non-positively curved manifolds into $M$ 
{\em parallel} if there
are a totally-geodesic embedding $\Phi :\tilde{\Sigma}_1\times
[a_1,a_2]\hookrightarrow \tilde{M}$ and Riemannian covering maps 
$p_i:\tilde{\Sigma}_1\times \{a_i\}\to \Sigma _i$ such that 
$\phi _i\circ p_i=\pi \circ\Phi \mid _{\tilde{\Sigma}_1\times \{a_i\}}$.
\end{dfn}

If we have two such immersions of $\Sigma $ into $M$ which are homotopic,
then there will be totally-geodesic submanifolds 
$Y$ and $Y^{\prime }$ in $\tilde M$ covering these
immersions and lying a bounded distance apart. 
Thus $Y$ and $Y^{\prime }$
are parallel, and hence so are the two immersions of $\Sigma $ into $M$.

Our aim is to prove the following result.

\begin{prop}
\label{finiteness} 
Let $M$ be a compact non-positively curved Riemannian
manifold with convex boundary. Then there is an upper bound to the number of
disjoint, closed, totally-geodesically embedded hypersurfaces in $M$ so that
no two of them are parallel.
\end{prop}

As discussed above, any such hypersurface is $\pi _1$-injective and if two
such hypersurfaces are homotopic, they must be parallel. Now in the
topological setting in dimension three, it is a standard result \cite{Hempel}
that, in any compact 3-manifold $M$, there is an upper bound to the number
of disjoint, embedded, $\pi _1$-injective closed surfaces in $M$ which are
pairwise non-parallel, where two surfaces $S$ and $S^{\prime }$ are parallel
if they together bound a submanifold homeomorphic to $S\times I.$ This upper
bound is called the Haken number of $M$. In higher dimensions, there is no
such result in the general topological setting, but there is an algebraic
analogue due to Dunwoody \cite{Dunwoody}, which discusses splittings of $PDn$%
-groups over $PD(n-1)$-subgroups. This implies that if one considers a
closed aspherical manifold $M$, there is an upper bound to the number of
disjoint closed aspherical $\pi _1$-injective embedded codimension-one
submanifolds in $M$ such that no two are homotopic. Clearly this result will
also apply to any aspherical compact manifold with boundary so long as the
boundary is also $\pi _1$-injective and aspherical. Now the hypotheses of
the above proposition imply that $M$ is aspherical and that its boundary is $%
\pi _1$-injective and aspherical. Thus one can prove this proposition in
dimension three by using the Haken number, and can prove it in any dimension
by using Dunwoody's result. However, we will give a direct geometric proof.

\medskip\noindent{\it Proof of Proposition \ref{finiteness}:} 
Let $\Sigma _1,\ldots ,\Sigma _n$ be such a
family of hypersurfaces. 
It is possible that some of these hypersurfaces are components of $\partial M$. 
Consider a component $N$ of 
$M\setminus \cup_{i=1}^n\Sigma _i$. 
Identify the universal cover $\tilde N$ with a component
of $\pi ^{-1}(N)\subset \tilde M$. 
We will use $\partial \tilde
N$ to denote the boundary of $\overline{\tilde N}$ as a manifold. 
Then $\partial\tilde N$ consists of a union of components
of $\partial \tilde M$ and of totally geodesic hypersurfaces in $\tilde M$
each of which covers one of the surfaces $\Sigma _i$. 
Also $\partial \tilde N$ consists of at least two components. 
For
otherwise, $N$ would have infinite diameter contradicting the compactness
of $M$. 

{\it \ }Consider the case that $\partial \tilde N$ has exactly two
components $Y_1$ and $Y_2$. At least one of them, say $Y_1$, covers a
hypersurface which we denote $\Sigma _1$. The distance function $d(\cdot ,Y_1)$ is bounded
on $Y_2$, and vice versa, because the subgroup of $\Gamma $ preserving $%
\tilde N$ and $Y_1$ also preserves $Y_2$ and acts cocompactly on $\overline{%
\tilde N}$. 
Thus $Y_1$ and $Y_2$ are parallel and $\overline{\tilde N}$ is isometric 
to a product $Y_1\times[-a,a]$. 
If also $Y_2$ covers a hypersurface $\Si_i$, then $\Si_i$ is parallel to 
$\Si_1$ and hence $\Si_i=\Si_1$ by our assumption. 
If $\overline{\tilde N}$ projects onto $M$, then $n=1$. 
Otherwise the image of $\overline{\tilde N}$ 
is a twisted interval bundle over the hypersurface $\Si'_1$ covered by $Y_1\times\{0\}$. 
In this case, we replace $\Si_1$ by $\Si'_1$. 
If $Y_2$ covers no hypersurface $\Si_i$ and hence $Y_2\subset\partial\tilde M$,
then we remove $N$ from $M$. 
In both cases,
this reduces by one the number of components of 
$M\setminus \cup_{i=1}^n\Sigma _i$ 
and does not alter the number of hypersurfaces $\Si_i$. 
By repeating these steps,  
we may 
assume that for all pieces $N$, the universal cover $\tilde N$ has at least
three boundary components. The pieces then have a certain minimal size:

\begin{lem}
\label{minimal size of components} Each component $N$ of $M\setminus \cup
_{i=1}^n\Sigma _i$ contains a point $p$ at distance at least $\rho _0$ from
the boundary $\partial N$, where $\rho _0$ is a positive constant only
depending on the lower sectional curvature bound of $M$.
\end{lem}

\medskip\noindent{\it Proof: } We re-scale so that the sectional curvature
of $M$ is bounded by $-1\leq K_M\leq 0$. Let $\tilde{p}\in \tilde{N}$ be a
point at maximal distance $\rho $ from $\partial \tilde{N}$. The ball $B$ of
radius $\rho $ centered at $\tilde{p}$ touches three components $Y_1,Y_2,Y_3$
of $\partial \tilde{N}$ in respective points $\tilde{p}_1,\tilde{p}_2,\tilde{%
p}_3$. Let $v_i$ be the unit vector in $\tilde{p}$ pointing in the direction
of $\tilde{p}_i$. Among the vectors $v_1,v_2,v_3$ at least two, say $v_1$
and $v_2$, enclose an angle $\angle (v_1,v_2)\leq {\frac 23}\pi $. Consider
the arc in the unit sphere in $T_{\tilde{p}}\tilde{M}$ joining $v_1$ and $%
v_2 $. It contains a vector $v$ such that the geodesic ray $r:[0,\infty )\to 
\tilde{M}$ emanating from $\tilde{p}$ in the direction of $v$ intersects
neither $Y_1$ nor $Y_2$. We assume without loss of generality that $\angle
(v,v_1)\leq {\frac 13}\pi $. The angles of the triangles $\tilde{p}\tilde{p}%
_1r(t)$ satisfy for all $t>0$: 
\[
\label{angle inequalities}\angle _{\tilde{p}}(\tilde{p}_1,r(t))\leq {\frac 13%
}\pi ,\qquad \angle _{\tilde{p}_1}(\tilde{p},r(t))\leq {\frac 12}\pi 
\]
Consider comparison triangles with the same side lengths in the hyperbolic
plane ${\Bbb H}^2$. By Toponogov's triangle comparison theorem~\cite{Karcher}%
, the angles in the comparison triangles are not greater than the
corresponding angles in the triangles $\tilde{p}\tilde{p}_1r(t)$. So they
satisfy analogous inequalities. Since $t$ may be arbitrarily large, we can
bound $\rho $ from below by a positive constant $\rho _0$, namely by the
finite sidelength of the triangle in ${\Bbb H}^2$ with angles $0,{\frac 13}%
\pi ,{\frac 12}\pi $ and one ideal vertex.
\qed

Denote by $c(n)$ the number of components of 
$M\setminus \cup_{i=1}^n\Sigma _i$. 
According to Lemma~\ref{minimal size of components}, there is a $2\rho _0$%
-net in $M$ with one point in each component of 
$M\setminus \cup_{i=1}^n\Sigma _i$. By compactness of $M$, $c(n)$ 
stabilizes as $n$ tends to infinity. 
More precisely,
it can be bounded above in terms of the lower curvature bound and the volume 
of $M$. 
If $c(n_1)=c(n_2)$ for $n_1<n_2$, 
then we can choose for each $n$ with $n_1<n\leq n_2$ 
a closed smooth path $\alpha _n$
which does not intersect $\Sigma _1,\dots ,\Sigma _{n-1}$ but intersects $%
\Sigma _n$ once transversally. 
Looking at the intersection numbers modulo 2
of the paths $\alpha _n$ with the surfaces $\Sigma _n$, we see that the 
$\alpha _n$ represent linearly independent homology classes in 
$H_1(M,{\Bbb Z}/2{\Bbb Z})$. 
Since $M$ is compact, we conclude that $n_2-n_1$ is bounded 
in terms of the topology of $M$. 
This completes the proof of Proposition \ref{finiteness}.
\qed

\section{Geometric decomposition along closed submanifolds}
\label{rigid geometry}

In this section, $M$ will always denote a compact, connected, non-positively
curved Riemannian manifold 
of dimension at least 3 which has convex boundary.
We investigate how the pattern of closed totally-geodesic 
flat hypersurfaces in $M$ 
is organized to yield a canonical geometric decomposition. 
In dimension three this is a geometric relization 
of the canonical topological decomposition 
due to Jaco, Shalen and Johannson. 
The decomposition of $M$ will be obtained by cutting along hypersurfaces 
of the following kind 
(see section \ref{decomposing}): 

\begin{dfn}
\label{isolated surface} 
A totally-geodesically immersed, closed, flat
hypersurface in $M$ is called {\bf isolated} if it does not intersect any
such hypersurface transversally.
\end{dfn}

Note that the definition also excludes
self-intersections. It is immediate that isolated closed flat hypersurfaces cover
embedded hypersurfaces and the images of two of them must coincide or be
disjoint.

Denote by $\Ga$ the fundamental group of $M$
thought of as a group of deck transformations acting on $\tilde M$. 
Closed flat hypersurfaces in $M$ are covered by $(n-1)$-flats in $\tilde M$ 
which are periodic in the sense of:

\begin{dfn}
A {\bf $\Gamma$-periodic flat} or {\bf $\Gamma $-flat} 
is a flat $F$ in $\tilde{M}$ such that the subgroup $\Gamma _F$ of $\Gamma$ 
preserving $F$ acts cocompactly on $F$. 
We call $F$ {\bf isolated} if it intersects no other 
$\Gamma$-flat transversally. 
\end{dfn}

Unless explicitely stated otherwise,
all flats considerd in this section
will be $(n-1)$-dimensional. 
A totally-geodesically immersed closed flat hypersurface 
is isolated if and only
if it is covered by isolated $\Gamma $-flats in $\tilde{M}$.

\subsection{Seifert fibred submanifolds}
\label{Seifert}

>From now on we will assume that $M$ is not closed and flat. 
We prove in this section that intersecting, totally-geodesically immersed,
closed, flat hypersurfaces in $M$ span a submanifold which is foliated by 
parallel closed flat submanifolds of codimension two. 
In dimension three, this foliation is a Seifert fibration 
by closed geodesics. Our arguments are
closely related to those in Casson's proof of the Torus Theorem 
in dimension three \cite{Casson}, 
but are simpler because of the curvature assumption which we are imposing
on the metric of $M$.

Let $A\subset\Ga$ be 
a free
abelian subgroup of rank $n-2$. 
Recall from section \ref{prelim} 
that the normaliser $N(A)$ of $A$ in $\Gamma$ 
acts cocmpactly on the set  
of minimal displacement $\mathop{\hbox{\rm MIN}}(A)$ and 
preserves its metric splitting (\ref{splittingofminsetofabeliabgroup}). 
The induced action of $A$ on $E$ is
cocompact, so $E$ is Euclidean space of dimension $n-2$.

Now let $H_A$ denote the closed convex hull 
of the union of all $A$-invariant 
$\Gamma $-flats. 
$H_A$ is $A$-invariant 
and hence has the form  
$H_A =Z\times\R\subseteq Y\times {\Bbb R}=\mathop{\hbox{\rm MIN}}(A)$ 
for a closed convex subset $Z$ of $Y$. 
Furthermore $H_A$ is preserved by $N(A)$
and Lemma~\ref{centraliseractscoco} 
implies that the action of $N(A)$ on $H_A$ is cocompact. 
Note that it is possible that $H_A$
consists of a single $A$-invariant $\Gamma $-flat and so has empty
interior.
In this case it will be convenient to write $\partial H_A=H_A$. 

\begin{lem}
\label{boundary of convex hull of flats generalised}
The boundary $\partial H_A$ is a disjoint union of 
$\Gamma$-flats.
\end{lem}

\medskip\noindent{\it Proof: } 
Each $A$-invariant $\Gamma $-flat
projects to a complete geodesic in $Y$.
Let $Z$ denote the closed convex hull of the
family ${\cal F}$ of all such geodesics, 
so that $Z$ is 
either a geodesic or a convex subset of $Y$ with non-empty interior 
whose boundary $\partial Z$ is a union of disjoint complete geodesics. 
Consequently, $\partial H_A$ is a disjoint union of $A$-invariant 
$(n-1)$-flats lying above $\partial Z$. 
According to Lemma~\ref{centraliseractscoco} the quotient
manifold $N(A) \backslash H_A $ is compact and therefore also its
closed subset $\partial (N(A) \backslash H_A )=N(A) \backslash
\partial H_A$. 
Hence the components of $\partial H_A$ are $\Gamma $-flats.
\qed

Next we consider how codimension-one flats in $\tilde M$ can meet $%
\mathop{\hbox{\rm MIN}}(A)$. 

\begin{lem}
\label{no transversal flat generalised}
Suppose that $\mathop{\hbox{\rm MIN}}(A)$ contains two non-parallel $A$%
-invariant $\Gamma $-flats. Then any codimension-one flat intersecting $\mathop{\hbox{\rm MIN}}(A)$ is
also $A$-invariant and so is completely contained in $\mathop{\hbox{\rm
MIN}}(A)$.
\end{lem}

\medskip\noindent{\it Proof of Lemma \ref{no transversal flat generalised}: } 
Let $F$ be a $(n-1)$-flat in $\tilde M$ 
and denote by $U\subset Y$ the image of $F\cap MIN(A)$ under the canonical 
projection $MIN(A)=E\times Y\to Y$. 
$U$ is a convex subset with the property 
that every geodesic segment $\si\subset Y$ 
which intersects $U$ in more than one point is contained in $U$. 
If $F$ intersects $\mathop{\hbox{\rm MIN}}(\gamma )$ 
and is not $A$-invariant then $F$ intersects some $A$-flat $E\times\{y\}$
transversally. 
Hence $U\subset Y$ has non-empty interior and therefore $U=Y$. 
Lemma~\ref{flats in products} implies that $Y$ is flat. 
By our assumptions, $Y$ contains two complete non-parallel geodesics.
Therefore $Y$ is isometric to Euclidean plane and $\tilde M$ 
is isometric to Euclidean $n$-space. 
This implies that $M$ is closed and flat which contradicts the assumption 
made at the beginning of this section. 
\qed

The above result is a key step in our argument, and a very similar result appears in
Casson's proof of the Torus Theorem in the 3-dimensional case \cite{Casson}. 
In Casson's argument, no
assumption is made about the metric on $M$. Instead of considering a totally
geodesic immersion of the torus in $M$, he considers a least area immersion.
This means that in the universal cover of $M$, he considers area minimizing
planes rather than flats. Any two such planes must be disjoint or intersect
transversely in a single line. Again this situation is very similar to that
in this paper, but double lines of area minimising planes need not be
geodesics. Call two of these double lines {\it weakly parallel} if there is
a non-trivial element of $\pi _1(M)$ which stabilises both of them. The
analogue of our lemma is his result that either $\pi _1(M)$ contains the
free abelian group of rank three, so that $M$ is closed and must admit a
flat metric, or that all the double lines are weakly parallel.

This result shows that under the hypotheses of 
\ref{no transversal flat generalised}, no 
$\Gamma $-flat can cross $\partial H_A$ transversally, as such a flat
would have to be $A$-invariant and so be completely contained in $%
H_A$. 
With Lemma~\ref{boundary of convex hull of flats generalised}, we obtain:

\begin{prop}
\label{boundary of convex hull of flats is isolated} 
If $\mathop{%
\hbox{\rm
MIN}}(A)$ contains two non-parallel $A$-invariant $\Gamma $%
-flats, then $H_A$ has non-empty interior and the boundary $\partial
H_A$ is a disjoint union of isolated $A$-invariant $\Gamma $-flats.
\end{prop}

The quotient of $H_A$ by $N(A)$ is a Seifert fibred manifold $S_A$
with fibres being closed flat manifolds of dimension $n-2$, 
and the fibres form a
totally geodesic foliation of $S_A$.
(The definition of Seifert fibered manifolds in arbitrary dimension 
is given at the end of section \ref{prelim}.) 
In dimension three, this is a foliation by closed geodesics.

\subsection{The decomposition}
\label{decomposing}

We continue to assume that $M$ is not a closed flat manifold.  
By Proposition~\ref{finiteness}, 
there are finitely many families of 
parallel isolated flat closed hypersurfaces in $M$. 
In order to avoid unnecessary flat pieces 
(which are topologically interval bundles over closed flat $(n-1)$-manifolds)
in the decomposition of $M$
obtained below, 
we choose in each family of parallel hypersurfaces 
a canonical one as follows: 
For a $\Ga$-flat $F\subset\tilde M$ the set of all $\Ga$-flats 
parallel to $F$ splits as a product $F\times I$ where $I$ is a closed 
connected subset of $\R$. 
Since $M$ is assumed not to be flat, 
$I$ is isometric to a compact interval $[-a,a]$.
\begin{dfn}
We call the $\Ga$-flat $F\times\{0\}$ and the 
immersed hypersurface which it coveres {\bf central}. 
We call the isolated $\Ga$-flat $F\subset\tilde M$ 
and the embedded hypersurface covered by it in $M$ {\bf preferred}
if either $F\subset\partial\tilde M$ or $F$ is central and 
$F\subset Int(\tilde M)$. 
\end{dfn}
A preferred $\Ga$-flat $F$ has the useful 
property that every $(n-1)$-flat $F'$
parallel to $F$ satisfies
$Stab_{\Ga}(F')  \subseteq  Stab_{\Ga}(F)$.
Accordingly, each isolated closed hypersurface can be homotoped 
to the unique preferred hypersurface parallel to it.  

We now consider the finite collection ${\cal F}$ of all preferred isolated 
hypersurfaces in $M$. 
They are disjoint, embedded and pairwise non-parallel. 
They decompose $M$ into finitely many pieces which are compact
non-flat manifolds with convex boundary. 
Let $N$ be a piece of the decomposition and denote its
fundamental group by $\Gamma ^{\prime }:=\pi _1(N)$. 
$N$ has the property that 
all its preferred isolated flat hypersurfaces are contained in the boundary
and hence every {\em isolated} closed flat hypersurface 
can be homotoped into the boundary. 
We have the following dichotomy:

\begin{itemize}
\item  All immersed closed flat hypersurfaces are isolated and can be 
homotoped into the boundary. 

\item $N$ contains non-isolated closed flat immersed hypersurfaces. 
\end{itemize}

This dichotomy corresponds to the two types of pieces occurring in the
topological decomposition theorem 
in the three-dimensional case, compare section 1. 
The pieces of the first kind are codimension-one atoroidal.
(See section \ref{prelim} for a definition; 
in dimension three this is equivalent to being atoroidal.)
Assume that $N$ is a piece of the second kind. 
Then $\tilde N$ contains two $\Ga$-flats $F_1$ and $F_2$ 
which intersect transversally in a $(n-2)$-flat $L$. 
$Stab_{\Ga'}(F_1)\cap Stab_{\Ga'}(F_2)$ acts cocompactly on $L$
by Lemma~\ref{cocompactintersections},
and it contains an abelian subgroup $A$ of finite index and rank $n-2$. 
According to \ref{boundary of convex hull of flats is isolated},
the corresponding Seifert fibered manifold $S_A$ has non-empty interior.
Each boundary component of $S_A$ is an isolated flat hypersurface 
and can hence be homotoped into $\partial N$. 
By the construction of $S_A$ it follows that $\partial S_A\subseteq\partial N$
and therefore $S_A=N$. 
Thus $N$ is a geometric Seifert piece. This concludes the
proof of the following result:

\begin{gdt}
Let $M$ be a compact connected non-positively curved manifold which has
convex boundary. Then either $M$ is closed and flat or the following holds.

Let ${\cal F}$ be the family of all preferred isolated
totally-geodesic closed flat co\-di\-men\-sion-one submanifolds of $M$.
Then ${\cal F}$ is a finite collection of disjoint, mutually non-parallel, 
embedded hypersurfaces 
and decomposes $M$ into compact manifolds with convex
boundary which are Seifert or atoroidal. The Seifert components are foliated
by codimension-two totally geodesic closed flat submanifolds 
and the foliation is
locally a Riemannian product foliation. Further any $\pi _1$-injective map
of a closed flat $(n-1)$-manifold into $M$\ can be homotoped to a totally
geodesic flat immersion, and any such immersion must lie in one of the
Seifert pieces or be parallel to a hypersurface of ${\cal F}$. 
\end{gdt}

\section{Splitting along submanifolds with boundary}

\label{general case}In this section we will state and prove our most general
result which corresponds to the full JSJ decomposition of a compact $3$%
-manifold with boundary.

If a Riemannian manifold $\Sigma $ has totally geodesic boundary, we will
abbreviate this to say that $\Sigma $ has TGB. 
A proper map into an $n$-manifold $M$ of a compact flat $(n-1)$%
-manifold with TGB is {\it essential} if it is $\pi _1$-injective and not
properly homotopic into the boundary of $M$.
We will say that $M$ is {\it simple} if it does not admit 
an essential map of a compact flat $(n-1)$-manifold with TGB. 

In order to prove our general decomposition theorem, we will consider a
compact connected non-positively curved manifold $M$ of dimension $n\geq 3$,
which has TGB. This assumption on the boundary means
that we can double $M$ along its boundary to obtain a closed connected
non-positively curved manifold $DM$ of dimension $n$. 
If $M$ is not flat, neither is $DM$ and 
we can apply our main
geometric decomposition theorem from the preceding section to obtain the
canonical decomposition of $DM$ by finitely many totally-geo\-de\-si\-c flat
closed submanifolds of codimension one. The fact that this splitting is
canonical means that it is invariant under the involution $\tau $ which
interchanges the two copies of $M$ in $DM$. Thus the intersection with $M$
of the canonical family of totally geodesic flat closed codimension-one
submanifolds of $DM$ yields the required canonical splitting of $M$. The
non-simple pieces of $M$ are obtained from the Seifert manifold pieces of $DM
$ by intersecting them with $M$. Thus these pieces of $M$ are either Seifert
manifolds themselves or they are ``half a Seifert manifold''. This second
case occurs when a Seifert piece $\Sigma $ of $DM$ is $\tau $-invariant, so
that the intersection of $\Sigma $ with $M$ consists of half of $\Sigma $.
The restriction of $\tau $ to $\Sigma $ is itself an isometry and it fixes $%
\Sigma \cap \partial M$ pointwise. Thus, for each component $\Omega $ of $%
\Sigma \cap \partial M$, this isometry of $\Sigma $ lifts to an isometry of
the universal cover of $\Sigma $ which fixes pointwise a copy $\Pi $ of the
universal cover of $\Omega $. Recall that the universal cover of $\Sigma $
is metrically a product $Z\times E$, where $Z$ is some 2-dimensional space
and $E$ is isometric to Euclidean space of dimension $n-2$. Also recall that 
$\Pi $ is part of a flat in the universal cover of $DM$ and hence is a flat
in $Z\times E$. It
follows that $\Pi $ is of the form $P\times Q$, where $P$ is some subset of $%
Z$ and $Q$ is some subset of $E.$ Hence $\Sigma \cap \partial M$ is either
vertical or horizontal in $\Sigma $, where (as in dimension three) a
codimension-one submanifold is vertical if it is a union of fibers of the
Seifert structure, and is horizontal if it is transverse to every fibre. In
the vertical case, $\Sigma \cap M$ is again a Seifert fiber space. In the
horizontal case, $\Sigma \cap M$ must be the product (or twisted product) of 
$\Sigma \cap \partial M$ with an interval.

\medskip
\noindent
{\bf General Geometric Decomposition Theorem.} {\em Let }$M$ {\em be a
compact connected non-positively curved manifold of dimension }$n\geq 3${\em %
, which has TGB. Then either the metric on }$M$ {\em is flat, or }$M${\em \ can be canonically
decomposed along finitely many totally-geo\-de\-si\-cal\-ly properly
embedded flat compact submanifolds of codimension one with TGB. The
resulting pieces are simple or Seifert fibered or are bundles with fibre an
interval over a compact }$(n-1)${\em -manifold with TGB. } {\em Further
any essential map into }$M${\em \ of a compact flat }$(n-1)${\em -manifold
with TGB\ can be properly homotoped to a totally geodesic flat immersion,
and any such immersion must lie in one of the non-simple pieces
or be parallel to one of the canonical surfaces.}

\noindent
Bernhard Leeb, Mathematisches Institut der Universit\"{a}t, 
Beringstr.\ 1, 53115 Bonn, Germany, 
leeb@rhein.iam.uni-bonn.de 

\noindent
Peter Scott,
Mathematics Department, 
University of Michigan, 
Ann Arbor, MI 48109, USA, 
pscott@umich.edu

\end{document}